\documentclass[11pt,openbib]{article} 
\usepackage{amsfonts,amsthm}
\usepackage{amssymb,amsmath}
\usepackage{epic}
\usepackage{graphicx}

\newcommand{\dee}{\mathrm{d}}

\newcommand{\setrule}{\, \mathop{\rule[-4pt]{.5pt}{13pt}\, }\nolimits}

\newcommand{\onehalf}{\mbox{$\frac{\scriptstyle 1}{\scriptstyle 2}\,$}} 
 
\newcommand{\ttfrac}[2]{\mbox{$\frac{{\scriptstyle #1}}{{\scriptstyle #2}}$}}
\allowdisplaybreaks
\title{\textbf{Comment on a theorem of Hamilton}}
\author{Richard Cushman\thanks{email: rcushman@ucalgary.ca} \\ Department of Mathematics and Statistics, \\
University of Calgary, Calgary, Alberta, Canada}
\date{6 January 2015}
\begin{document}
\maketitle
\noindent This paper gives a slight refinement of a theorem of Hamilton \cite{hamilton}, which shows that the velocity of a Keplerian motion moves on a circle. \medskip 

\noindent The motion of a particle of mass $m$ under an attractive central force with potential $U(|x|) = 
-k \,  \frac{1}{|{\bf x}|}$ 
is governed by Newton's equations
\begin{equation}
m \frac{{\dee }^2 {\bf x}}{{\dee t}^2} = -\mathrm{grad}\, U(|{\bf x}|) = -k \, \frac{{\bf x}}{|{\bf x}|^3}.
\label{eq-one}
\end{equation}
Here ${\bf x} \in {\bf R}^3 \setminus \{ (0,0,0) \}$ with $|{\bf x}|$ being the length of ${\bf x}$ using the Euclidean 
inner product $\langle \, \, , \, \, \rangle $.\medskip

\noindent In \cite{hamilton} Hamilton proved \medskip 

\noindent \textbf{Theorem}. The velocity vector ${\bf v} = \frac{\dee {\bf x}}{\dee t}$ of the particle moves 
on a circle $\mathcal{C}$, which uniquely determines its Keplerian orbit. \medskip 

\noindent Following Milnor \cite{milnor}, see also Anosov \cite{anosov}, we recall the proof of Hamilton's theorem. \medskip 

\noindent \textbf{Proof}. Let ${\bf J} = {\bf x} \times m{\bf v}$ be the angular momentum of the particle, which we 
assume is nonzero. Because 
\begin{displaymath}
\frac{\dee {\bf J}}{\dee t} = \frac{\dee {\bf x}}{\dee t} \times m{\bf v} + {\bf x} \times m\frac{\dee {\bf v}}{\dee t} 
={\bf v} \times m {\bf v} - {\bf x} \times \big( -\frac{k}{|{\bf x}|^3}{\bf x} \big) = {\bf 0}, 
\end{displaymath}
it follows that ${\bf J}$ is a constant of motion. Since ${\bf J} \ne {\bf 0}$, the motion ${\bf x}(t)$ and the 
velocity ${\bf v}(t)$ of the particle are linearly independent vectors, which lie in, and thus span, a plane $\Pi $ perpendicular to ${\bf J}$. \medskip 

\noindent Introduce coordinates on ${\bf R}^3$ so that ${\bf J} = (0,0,j)$, where $j = |{\bf J}| > 0$. Then 
$\Pi =\{ {\bf x} = (x,y,0) \in {\bf R}^3 \, \setrule \, (x,y) \in {\bf R}^3 \}$. Using polar coordinates $(r, \theta )$ on 
$\Pi $ so ${\bf x} = (r \cos \theta , r \sin \theta , 0)$, we find that 
\begin{equation}
j = x \frac{\dee y}{\dee t} - y \frac{\dee x}{\dee t} = r^2 \, \frac{\dee \theta }{\dee t}. 
\label{eq-two}
\end{equation}
From (\ref{eq-two}) and the fact that $r >0$ and $j >0$, it follows that $\frac{\dee \theta }{\dee t} >0$. Therefore 
we can reparametrize the curves $t \mapsto {\bf x}(t)$ and $t \mapsto {\bf v}(t)$ using $\theta $ instead of $t$. This reparametrization preserves the original positive orientation of both curves, given by increasing $t$. Now 
write Newton's equations (\ref{eq-one}) as 
\begin{displaymath}
\frac{\dee {\bf v}}{\dee t} = -\frac{k}{m r^2}\, ( \cos \theta , \sin \theta , 0).
\end{displaymath}
Dividing by $\frac{\dee \theta }{\dee t}$ and using (\ref{eq-two}) gives 
\begin{displaymath}
\frac{\dee {\bf v}}{\dee \theta } = \frac{\frac{\dee {\bf v}}{\dee t}}{\frac{\dee \theta }{\dee t}} = 
- \big( k/mr^2 \mbox{\LARGE $/$} j/r^2 \big) \, (\cos \theta , \sin \theta , 0) = - R \, (\cos \theta , \sin \theta , 0), 
\end{displaymath}
where $R = k/jm$. Integrating the above equation we get 
\begin{equation}
{\bf v}(\theta ) = R \, (-\sin \theta , \cos \theta , 0) + {\bf c}, 
\label{eq-three}
\end{equation}
where ${\bf c} = (c_1, c_2, 0)$. Thus the velocity vector ${\bf v}$ of the particle moves in a positive sense, namely, 
with increasing $\theta $, on a circle $\mathcal{C}$ in the plane $\Pi$ with center at ${\bf c}$ and radius $R$. 
\medskip  

\noindent Let $e = c/R$, where $c = |{\bf c}| \ge 0$. Choose coordinates on 
$\Pi $ so that ${\bf c} = (0,c,0)$. We may rewrite (\ref{eq-three}) as 
\begin{equation}
{\bf v}(\theta ) = \big( -R \sin \theta , R(e+ \cos \theta ) , 0  \big) . 
\label{eq-four}
\end{equation}
Consequently, 
\begin{align}
j & = \langle {\bf J}, (0,0,1) \rangle = \langle {\bf x}(\theta ) \times m {\bf v}(\theta ), (0,0,1) \rangle \notag \\
& = \big( r(\theta ) \cos \theta \big) mR (e +\cos \theta ) - \big( r(\theta ) \sin \theta \big) mR (-\sin \theta ) \notag \\
& \hspace{1in}\mbox{using ${\bf x}(\theta ) = \big( r(\theta ) \cos \theta , r(\theta )\sin \theta , 0)$ and 
(\ref{eq-four}) } \notag \\
& = m r(\theta ) R (1 + e \cos \theta ) . \notag 
\end{align}
So the orbit on $\Pi $ traced out by the motion $\theta \mapsto {\bf x}(\theta )$ of the particle satsifies 
\begin{equation}
r = r(\theta ) = \Lambda (1+e \cos \theta )^{-1}. 
\label{eq-five}
\end{equation}
This is the equation of a conic section of eccentricity $e \ge 0$ with focus at $O = (0,0,0)$. Here 
$\Lambda = j/mR = j^2/k$. \hfill $\square $ \medskip 

\noindent In the case of elliptical or circular motion $0 \le e < 1$ the angle $\theta $ increases by $2\pi $, while 
the particle traces out the ellipse or circle, respectively. Hence the velocity vector traces out all of the velocity circle 
$\mathcal{C}$. \medskip 

\noindent From now on we only 
look at the case of hyperbolic Keplerian motion (\ref{eq-five}) where 
$e > 1$. So $|\theta | < {\theta }_0 = {\cos }^{-1}\big( -e^{-1} \big) = \pi - {\theta }_{\ast }$, where ${\theta }_{\ast } = 
{\cos }^{-1}e^{-1}$. In this case we will show that the velocity vector traces out a closed arc $\mathcal{A}$ of the 
circle $\mathcal{C}$ in a positive sense. Conservation of energy places the following constraint 
\begin{equation}
|{\bf v}|^2 = \langle {\bf v}, {\bf v} \rangle = \frac{2h}{m} +\frac{k}{m^2} \, \frac{1}{|{\bf x}|} > \frac{2h}{m} 
\label{eq-six}
\end{equation}
on the length squared of the velocity of the particle. In other words, the velocity vector ${\bf v}$ lies outside of the closed $2$-disk $\mathcal{E}$ in $\Pi $ with center at $O$ and 
radius $\sqrt{\frac{2h}{m}}$. It lies on the velocity circle $\mathcal{C}$ and on the energy circle $\partial \mathcal{E}$, given by $|{\bf v}| = \sqrt{\frac{2h}{m}}$, if and only if $0 = \frac{1}{r(\theta )} = {\Lambda }^{-1}(1+e\cos \theta )$, that is, 
if and only if $\theta = \pm {\theta }_0 = \pm (\pi - {\theta }_{\ast})$. Thus the closed arc $\mathcal{A}$ has 
end points ${\bf v}(\pm {\theta }_0)-{\bf c}$. Using (\ref{eq-four}) we see that the 
corresponding velocity vector on $\mathcal{C}$ is 
\begin{align}
{\bf v}(\pm {\theta }_0) - {\bf c} & = \big( -R \sin (\pm {\theta }_0), R \cos (\pm {\theta }_0), 0 \big) = 
\big( \mp R \sin {\theta }_{\ast} , -R \cos {\theta }_{\ast }, 0 \big) \notag \\
&\hspace{-.6in} = \left\{ \begin{array}{rl} 
\big( R \cos (\ttfrac{3}{2} \pi  - {\theta }_{\ast }), R \sin (\ttfrac{3}{2} \pi - {\theta }_{\ast }), 0 \big), & \mbox{when $+$ holds} \\
\big( R \cos (-(\ttfrac{1}{2} \pi - {\theta }_{\ast })), R \sin (-(\ttfrac{1}{2} \pi  - {\theta }_{\ast })), 0 \big), & \mbox{when $-$ holds.} 
\end{array} \right. 
\label{eq-seven}
\end{align}
Now ${\bf v}(\pm {\theta }_0)$ is the asymptotic velocity of the outgoing motion of the particle when the $+$ sign is taken, and the asymptotic velocity of the incoming motion when the $-$ sign is taken. Note that ${\bf v}({\theta }_0)-{\bf c}$ lies in the intersection of the half planes $\{ x < 0 \} $ and $\{ y < 0 \} $ of the $2$-plane $\Pi $; while ${\bf v}(-{\theta }_0)$ lies in the intersection of the half planes $\{ x > 0 \} $ and $\{ y < 0 \} $ and is symmetric in the $y$-axis to 
${\bf v}({\theta }_0)$. Thus the velocity of the particle moves along the arc $\mathcal{A}$ of 
$\mathcal{C}$ in the \emph{positive} sense (with $\theta $ \emph{increasing}) from ${\bf v}(-{\theta }_0)-{\bf c}$ to 
${\bf v}({\theta }_0) - {\bf c}$. Let $\Theta $ be the positive angle swept out by a counterclockwise rotation about ${\bf c}$, which sends ${\bf v}(-{\theta }_0)- {\bf c}$ to ${\bf v}({\theta }_0)-{\bf c}$. From (\ref{eq-seven}) it follows that 
$\Theta = 2(\pi - {\theta }_{\ast })$. The directions of the asymptotic motion of the particle corresponding to 
${\bf v}(\pm {\theta }_0)$ are 
\begin{equation}
{\bf d}_{\pm {\theta }_0} = \frac{{\bf x}(\pm {\theta }_0)}{|{\bf x}(\pm {\theta }_0)|} = 
\big( \cos (\pm {\theta }_0), \sin (\pm {\theta }_0), 0 \big) = \big( -\cos {\theta }_{\ast }, \pm \sin {\theta }_{\ast }, 0 \big). 
\label{eq-eight}
\end{equation}
Here ${\bf d}_{{\theta }_0}$ is the asymptotic direction of the outgoing motion of the particle; while ${\bf d}_{-{\theta }_0}$ 
is the asymptotic direction of incoming motion. By definition, the scattering angle $\Psi $ of the hyperbolic 
motion of the particle is the positive angle swept out by a counterclockwise rotation about the center $C =(ae,0,0)$ of the 
hyperbola which sends ${\bf d}_{-{\theta }_0}$ to ${\bf d}_{{\theta }_0}$.  \medskip 

\noindent \textbf{Claim}. In the case of hyperbolic Keplerian motion of energy $h$ and angular momentum of 
magnitude $j$, the angle $\Theta $ determined by the positive arc $\mathcal{A}$ of the velocity circle 
$\mathcal{C}$ is equal to the scattering angle $\Psi $.  \medskip 

\noindent \textbf{Proof}. By definition $\onehalf \Psi $ is the angle swept out by a counterclockwise rotation about the center $C$ of the hyperbola from the $x$-axis of $\Pi $ to the outgoing asymptotic direction ${\bf d}_{{\theta }_0}$ of the hyperbola. By construction $\onehalf \Psi = {\theta }_0$. So 
$\onehalf \Psi = {\theta }_0 = \pi - {\theta }_{\ast } = \onehalf \Theta $. Explicitly, we have 
$\onehalf \Theta = \pi - {\tan }^{-1}(\sqrt{e^2-1}) = \pi -{\tan }^{-1}\big( \frac{j}{k}\sqrt{2hm} \big) $. 
To see this last equality, substitute (\ref{eq-four}) and (\ref{eq-five}) into the conservation of energy equation 
(\ref{eq-six}). We get 
\begin{align}
h & = \onehalf mR^2 \big( {\sin }^2\theta +(e+\cos \theta )^2 \big) - \frac{k}{m\Lambda} (1+e \cos \theta ) \notag \\
& = \onehalf mR^2 (1 +e^2 + 2e\cos \theta ) - mR^2 (1+e\cos \theta ) \notag \\
& = \onehalf mR^2(e^2-1) = \onehalf \frac{k^2}{mj^2}(e^2-1). \tag*{$\square $}
\end{align}

\end{document}